\documentclass[notitlepage,12pt]{article}
\usepackage{amssymb,amsmath,geometry
}

\catcode`\@=11
\@addtoreset{equation}{section}

\catcode`\@=12

\allowdisplaybreaks

\usepackage{color}

\def\R{\mathbb{R}}

\def\f{\varphi}
\def\iirn{\iint\limits_{\R^n\!\times\R^n}}
\def\iirnp{\iint\limits_{\R^n_+\!\times\R^n_+}}

\def\irn{\int\limits_{\R^n}}

\def\irnplus{\int\limits_{\R^n_+}}

\def\eps{\varepsilon}

\def\Ds{\left(-\Delta\right)^{\!s}\!}  

\def\sstar{{2^*_s}}

\def\proof{\noindent{\textbf{Proof. }}}
\def\QED{\hfill {$\square$}\goodbreak \medskip}

\newtheorem{Theorem}{Theorem}[section]
\newtheorem{Lemma}[Theorem]{Lemma}
\newtheorem{Proposition}[Theorem]{Proposition}

\begin{document}

\title 
{Fractional Hardy-Sobolev inequalities on half spaces}

\author{Roberta Musina\footnote{Dipartimento di Matematica ed Informatica, Universit\`a di Udine,
via delle Scienze, 206 -- 33100 Udine, Italy. Email: {roberta.musina@uniud.it}. 
{Partially supported by Miur-PRIN 2015KB9WPT.}}~ and
Alexander I. Nazarov\footnote{
St.Petersburg Department of Steklov Institute, Fontanka 27, St.Petersburg, 191023, Russia, 
and St.Petersburg State University, 
Universitetskii pr. 28, St.Petersburg, 198504, Russia. E-mail: al.il.nazarov@gmail.com.
Supported by  RFBR grant 17-01-00678.
}
}

\date{}

\maketitle

\begin{abstract}
{\small 
We investigate the existence of extremals for Hardy--Sobolev inequalities involving the 
Dirichlet fractional Laplacian $\Ds$ of order
$s\in(0,1)$ on half-spaces.}

\medskip

\medskip

\noindent
\textbf{2010 Mathematics Subject Classfication:} 47A63; 35A23.
\end{abstract}

\normalsize

\bigskip

\section{Introduction}
\label{S:Introduction}

We study Hardy-Sobolev type inequalities for the restricted Dirichlet fractional Laplacian $\Ds$
acting on functions that vanish outside an half-space, for instance outside 
$$\R^n_+=\{~\!x=(x_1,x')\in\R\times\R^{n-1}~|~~ x_1>0~\!\}.$$
We always assume $s\in(0,1)$, $n>2s$ and we put
$$
\sstar:=\frac{2n}{n-2s}~\!.
$$
We recall that the operator $\Ds$ is defined by 
$$
{\mathcal F}\big[\Ds u\big] = |\xi|^{2s}{\mathcal F}[u]~,\quad
u\in {\cal C}^\infty_0(\R^n),
$$
where ${\mathcal F}$ is the Fourier transform
$\displaystyle{{\mathcal F}[u](\xi)= (2\pi)^{-\frac{n}{2}}
\int_{\R^n} e^{-i~\!\!\xi\cdot x}u(x)~\!dx}$. The corresponding quadratic form is given by
$$
\langle\Ds u,u\rangle=\irn|\xi|^{2s}|{\mathcal F}[u]|^2~\!d\xi.
$$
Motivated by applications to variational fractional
equations on half-spaces, in the present paper we study the inequality 
\begin{equation}
\label{eq:general}
\langle\Ds u,u\rangle\ge \lambda \irnplus x_1^{-2s}|u|^2~\!dx+\mathcal S_s^{\lambda,p}(\R^n_+)
\Big(\irnplus x_1^{-pb}|u|^p~\!dx\Big)^{\frac 2p},~~u\in {\cal C}^\infty_0(\R^n_+)
\end{equation}
under the following hypotheses on the data:
\begin{subequations}
\begin{gather}
\label{eq:splambda}
        2<p\le\sstar~,\quad \lambda< \mathcal H_s:=\displaystyle
{\frac{1}{\pi}\Gamma\big(s+\frac12\big)^2}~\!\\
        \frac{b}{n}=\frac{1}{p}-\frac{1}{\sstar}~\!.
        \label{eq:assumption}
\end{gather}
\end{subequations}
The bounds on the exponent $p$ are due to Sobolev embeddings; 
the relation (\ref{eq:assumption}) is a necessary condition to have
of (\ref{eq:general}) for some constant $\mathcal S_s^{\lambda,p}(\R^n_+)>0$, use a rescaling argument.

 Actually the assumptions (\ref{eq:splambda}--\ref{eq:assumption}) are sufficient to have that 
(\ref{eq:general}) holds with a positive best constant $\mathcal S_s^{\lambda,p}(\R^n_+)$. Here is the argument.

First, notice that for $p=\sstar$, that implies $b=0$, we have
\begin{equation}
\label{eq:Sobolev}
\mathcal S_s:=\inf_{\scriptstyle u\in \mathcal C^\infty_0(\R^n)\atop\scriptstyle  u\ne 0}
\frac{\langle \Ds u,u\rangle}
{\|u\|^2_\sstar}=
\inf_{\scriptstyle u\in \mathcal C^\infty_0(\R^n_+)\atop\scriptstyle  u\ne 0}
\frac{\langle \Ds u,u\rangle}
{\|u\|^2_\sstar}= \mathcal S_s^{0,\sstar}(\R^n_+)
\end{equation}
because of the action of translations and dilations in $\R^n$.   The explicit value of the Sobolev constant
$\mathcal S_s$ has been computed in \cite{CoTa}.

Next, recall the 
Hardy-type inequality with cylindrical weights proved by Bogdan and Dyda  in \cite{BD}. It turns out that 
\begin{equation}
\label{eq:Hardy}
\langle \Ds u,u\rangle\ge \mathcal H_s \irnplus x_1^{-2s}u^2~\!dx\qquad\text{for any  $u\in\mathcal C^\infty_0(\R^n_+)$,}
\end{equation}
with a sharp constant in the right hand side. Thus $\mathcal S_s^{\lambda,\sstar}(\R^n_+)>0$ for any
$\lambda<\mathcal H_s$.

If $p\in(2,\sstar)$ and (\ref{eq:splambda}--\ref{eq:assumption}) are satisfied,
the existence of a positive constant
$\mathcal S_s^{\lambda,p}(\R^n_+)$  such that (\ref{eq:general}) holds  is easily proved via H\"older interpolation
between  the Sobolev and the cylindrical Hardy inequalities.

We now set up an appropriate functional setting to study
the existence of extremals for $\mathcal S_s^{\lambda,p}(\R^n_+)$. The quadratic form $\langle \Ds u,u\rangle$ 
induces an Hilbertian structure on the space
$$
{\mathcal D}^s(\R^n)\!=\!\{u\in {L^\sstar(\R^n)}~|~\langle \Ds u,u\rangle<\infty~\},
$$
and ${\mathcal D}^s(\R^n) \hookrightarrow L^\sstar(\R^n)$ with a continuous embedding by the Sobolev inequality.
Clearly ${\mathcal D}^s(\R^n)\cap L^2(\R^n)$ is the standard Sobolev space
$H^s(\R^n)$, see \cite{Tr} for basic results about $H^s$-spaces. In particular
${\mathcal D}^s(\R^n)\supsetneq H^s(\R^n)$ and ${\mathcal D}^s(\R^n)\subset H^s_{\rm loc}(\R^n)$,
that means $\f u\in H^s(\R^n)$ for  $\f\in\mathcal C^\infty_0(\R^n)$ and $u\in {\mathcal D}^s(\R^n)$.
Therefore, 
$\mathcal C^\infty_0(\R^n)$ is dense in ${\mathcal D}^s(\R^n)$ and  the {Rellich-Kondrashov Theorem} holds,
that is,
${\mathcal D}^s(\R^n)$ is compactly embedded into 
$L^q_{\rm loc}(\R^n)$ for any $q<\sstar$.

Next, let $\widetilde{\mathcal D}^s(\R^n_+)$ be the closure of $\mathcal C^\infty_0(\R^n_+)$ in ${\mathcal D}^s(\R^n)$.
We have
\begin{gather}
\nonumber
\widetilde{\mathcal D}^s(\R^n_+)=\{u\in \mathcal D^s(\R^n)~|~u\equiv 0~~\text{on} ~~\R^n_-:=
\R^n\setminus\overline\R^n_+~\}~\!,\\
\label{eq:problem}
\mathcal S_s^{\lambda,p}(\R^n_+)
=\inf_{\scriptstyle u\in \widetilde{\mathcal D}^s(\R^n_+)\atop\scriptstyle  u\ne 0}
\frac{\langle \Ds u,u\rangle-\lambda \|x_1^{-s}u\|^2_2}
{\|x_1^{-b}u\|^2_p}~\!.
\end{gather}
The minimization problem in (\ref{eq:problem}) is noncompact, due to the action
of dilations in $\R^n$. If $n\ge 2$ also translations in $\R^{n-1}$ might generate
noncompact minimizing sequences. Our first result about the existence of minimizers for
$\mathcal S_s^{\lambda,p}(\R^n_+)$ concerns the case $p<\sstar$.

\begin{Theorem}
\label{T:subcritical}
Let $p\in(2,\sstar)$, $-\infty<\lambda<\mathcal H_s$ and $b\in (0,s)$ as in (\ref{eq:assumption}).
Then the infimum $\mathcal S_s^{\lambda,p}(\R^n_+)$ is achieved in
$\widetilde{\mathcal D}^s(\R^n_+)$.
\end{Theorem}

In the critical case $p=\sstar$,  the noncompact group translations in the $x_1$-variable
produces severe lack of compactness phenomena. Take for instance $\lambda=0$. By the results in
\cite{CoTa} we have that, up to dilations, translations and multiplications, the Sobolev constant
$\mathcal S_s$ is attained on $\mathcal D^s(\R^n)$ only by the function
$
U_{\!s}(x)=\displaystyle{\left(1+|x|^2\right)^{\frac{2s-n}{2}}}$.
Therefore the infimum $\mathcal S_s^{0,\sstar}(\R^n_+)=\mathcal S_s$ is not achieved on 
$\widetilde{\mathcal D}^s(\R^n_+)$. In Section \ref{S:Dirichlet} we 
prove the next theorems.

\begin{Theorem}
\label{T:positive}
For $\lambda< \mathcal H_s$ the following facts hold.
\begin{itemize}
\item[$i)$] $\mathcal S_s^{\lambda,\sstar}(\R^n_+)\le  \mathcal S_s$;
\item[$ii)$] If $-\infty<\lambda~\!\le~\! 0$~ then $\mathcal S_s^{\lambda,\sstar}(\R^n_+)= \mathcal S_s$
and $\mathcal S_s^{\lambda,\sstar}(\R^n_+)$ is not achieved;
\item[$iii)$] If ~~~$0~\!<~\!\lambda<\mathcal H_s$ and $n\ge 4s$ then
$\mathcal S_s^{\lambda,\sstar}(\R^n_+)< S_{\!s}$.
\end{itemize}
\end{Theorem}

\begin{Theorem}
\label{T:Dirichlet} 
Assume $0<\lambda<{\mathcal H}_s$.  If $\mathcal S_s^{\lambda,\sstar}(\R^n_+)< S_{\!s}$ then 
$\mathcal S_s^{\lambda,\sstar}(\R^n_+)$
is achieved. In particular, if $n\ge 4s$ then $\mathcal S_s^{\lambda,\sstar}(\R^n_+)$ is achieved.
\end{Theorem}

Notice that $\mathcal S_s^{\lambda,\sstar}(\R^n_+)$ is always achieved if $n\ge 4$, while
the cases
$$
n=1~~\text{and~ $\frac14<s<\frac12$}~,\quad n=2~~\text{and ~$\frac12<s<1$}~,
\quad n=3~~\text{and ~$\frac34<s<1$}
$$
are not covered by Theorem \ref{T:Dirichlet}.

\medskip

All the proofs can be found in the next section. 
Our arguments to get the existence of minimizers are simple and self-contained.
We construct an ad hoc bounded minimizing sequence that can neither concentrate at the origin nor vanish.
In the locally compact case (see Theorem \ref{T:subcritical}) the existence of a minimizer is readily
obtained.
In the critical case, concentration at points $x\in\R^n_+$ is excluded by the assumption 
$\mathcal S_s^{\lambda,\sstar}(\R^n_+)< S_{\!s}$, and the existence result in Theorem
\ref{T:Dirichlet} follows.

Thanks to formula (\ref{eq:F}) below, an alternative proof can be 
obtained by adapting the arguments in the recent paper \cite{Fch}.

We conclude
the paper with few additional remaks and open problems. In particular, in Section
\ref{S:HSM} we conjecture that  Theorem \ref{T:Dirichlet} is sharp, that
is, $\mathcal S_s^{\lambda,\sstar}(\R^n_+)$ is not attained if $2s<n<4s$.

\bigskip
\noindent
{\bf Notation.}
{\small $\Omega\subset \R^n$ is a domain, and $\Omega^\mathsf{c}=\R^n\setminus\Omega$ is its complement.

For $q\in[1,\infty]$ we denote by $\|~\!\!\cdot~\!\!\|_{q,\Omega}$ the norm in 
$L^q(\Omega)$. If $\Omega=\R^n$ we simply write $\|~\!\!\cdot~\!\!\|_{q}$. 

Let $u\in \mathcal C^\infty_0(\R^n)$ and $s\in(0,1)$. It is well known that
\begin{equation}
\label{eq:Cns}
\langle \Ds u,u\rangle=\frac{C_{n,s}}{2}\iirn\frac{(u(x)-u(y))^2}{|x-y|^{n+2s}}~dxdy~,
\quad C_{n,s}=\displaystyle{\frac{s2^{2s}\Gamma\big(\frac{n}{2}+s\big)}{\pi^{\frac{n}{2}}\Gamma\big(1-s\big)}}~\!.
\end{equation}
By density, (\ref{eq:Cns}) holds for any $u\in \mathcal D^s(\R^n)$. Next, for $\lambda<\mathcal H_s$ we put
$$
\mathcal E^\lambda_s(u)=\langle \Ds u-\lambda x_1^{-2s} u,u\rangle=\langle\Ds u, u\rangle-\lambda\int_{\R^n_+} x_1^{-2s}|u|^2~\!dx~,
$$
that is the square of an equivalent norm in $\widetilde{\mathcal D}^s(\R^n_+)$ by the Hardy inequality (\ref{eq:Hardy}).

Through the paper, all constants depending only on $n$ and $s$ are denoted by $c$. To indicate that a constant depends on other quantities 
we list them in parentheses: $c(\dots)$.}

\section{Proofs}
\label{S:Dirichlet}

We start with a technical result that is essentially known, see for instance \cite{MNmp}. We 
provide its proof for the convenience of the reader.

\begin{Lemma}
\label{L:technicalD}
Let $u\in \widetilde{\mathcal D}^s(\R^n_+)$, $\f\in\mathcal C^\infty_0(\R^n)$ and let $\Omega
\subset\R^n$ be a bounded
 domain containing the support of $\f$. Then $\f u\in \widetilde{\mathcal D}^s(\R^n_+)$ and
 $$
\big|\langle \Ds \f u,\f u\rangle-\langle \Ds u,\f^2 u\rangle\big|\le  c(\f,\Omega)\langle \Ds u,u\rangle^{\frac 12}\cdot\|u\|_{2,\Omega}.
$$
\end{Lemma}

\proof
The first statement is evident. Further, we estimate 
$$
\Psi_\f(x,y):=\frac{(\f(x)-\f(y))^2}{|x-y|^{n+2s}}\le c(\f)\Big(\frac{\chi_{\{|x-y|<1\}}}{|x-y|^{n-2(1-s)}}+\frac{\chi_{\{|x-y|>1\}}}{|x-y|^{n+2s}}\Big)$$ 
to obtain  
\begin{equation}
\label{eq:Psi_estimate}
\irn\Psi_{\!\f}(x,y)~\!dy\le c(\f)
\end{equation}
for all $x\in\R^n$. Taking (\ref{eq:Cns}) into account, by direct computation one finds 
$$
\langle \Ds \f u,\f u\rangle-\langle\Ds u,\f^2 u\rangle=
c\iirn{u}(x){u}(y)\Psi_\f(x,y)~\!dxdy=:B_\f~\!.
$$
Since the support of $\Psi_\f$ is contained in $(\Omega\times\R^n)\cup(\R^n\times\Omega)$, we have
\begin{multline*}
c~\!\big|B_\f\big|\le \iint\limits_{\Omega\!\times\Omega} |u(x)u(y)|\Psi_\f(x,y)~\!dxdy+
\int\limits_{\Omega}|u(x)|\Big(\int\limits_{\Omega^\mathsf{c}} \frac{|u(y)||\f(x)|^2}{|x-y|^{n+2s}}dy
\Big)dx\\
\stackrel{*}{\le} \iint\limits_{\Omega\!\times\Omega} |u(x)u(y)|\Psi_\f(x,y)~\!dxdy
+\int\limits_{\Omega} |u(x)\f(x)|^2\Big(\int\limits_{\Omega^\mathsf{c}}\frac{dy}{|x-y|^{n+2s}}\Big)dx\\
+\|\f\|_\infty\iint\limits_{\Omega\!\times \Omega^\mathsf{c}} \frac{|u(x)-u(y)|}{|x-y|^{\frac{n+2s}{2}}}~\!\frac{|u(x)\f(x)|}{|x-y|^{\frac{n+2s}{2}}}~\!dxdy
=:I_1+I_2+\|\f\|_\infty I_3
\end{multline*}
(in ($*$) we use the triangle inequality). 
By the Cauchy-Bunyakovsky-Schwarz inequality and (\ref{eq:Psi_estimate}) we obtain
\begin{multline*}
I_1=  \iint\limits_{\Omega\!\times\Omega}\big||u(x)|^2\Psi_\f(x,y)\big|^\frac12~\!
\big||u(y)|^2\Psi_\f(x,y)\big|^\frac12~\!dxdy\\
\le  \iint\limits_{\Omega\!\times\Omega}|u(x)|^2\Psi_\f(x,y)~\!dxdy
\le  c(\f)\int\limits_{\Omega}|u(x)|^2~\!dx\le  c(\f,\Omega)\langle \Ds u,u\rangle^{\frac 12}\cdot\|u\|_{2,\Omega}.
\end{multline*}
Since $\text{supp}(\f)$ is compactly contained in $\Omega$ we clearly have 
$\displaystyle{I_2\le c(\f,\Omega) \int\limits_{\Omega}|u(x)|^2dx}$.
To handle $I_3$ we use  the Cauchy-Bunyakovsky-Schwarz inequality in $\Omega\times \Omega^\mathsf{c}$ 
and the above estimate on $I_2$ to get
$$
\displaystyle{I_3^2\le c ~\! \langle \Ds u,u\rangle I_2 
\le c(\f,\Omega) \langle \Ds u,u\rangle \int\limits_{\Omega}|u(x)|^2dx}.
$$
The proof is complete.
\QED

\bigskip
\noindent{\bf Proof of Theorem \ref{T:subcritical}.}
We follow the outline of the proof of Theorem 0.1 in \cite{GM}. 
Thanks to  a standard convexity argument, we only need to construct a minimizing sequence
that weakly converges in $\widetilde{\mathcal D}^s(\R^n_+)$ to a nontrivial limit. 
For future convenience we notice that the assumption $p<\sstar$
is only used in the last line of the proof. 

In order to simplify notations we put 
$${\mathcal S_\lambda}=\mathcal S^{\lambda,p}_s(\R^n_+)=
\inf_{\scriptstyle u\in \widetilde{\mathcal D}^s(\R^n_+)\atop\scriptstyle  u\ne 0}
\frac{\mathcal E_s^\lambda(u)}{\|x_1^{-b}u\|^2_p}.
$$
We assume that $n\ge 2$. The proof for $n=1$ is similar, and simpler; only notation
has to be adapted. 

For $\rho>0$ and ${z}\in \R^{n-1}$ we denote by $B'_\rho({z})$ the  $(n-1)$-dimensional ball
$$
B'_\rho({z})=\{x'\in\R^{n-1}~|~|x'-{z}|<\rho~\}~\!.
$$
Choose a finite number of points  
$x'_1,\cdots x'_\tau\in \R^{n-1}$ such that
\begin{equation}
\label{eq:Rcover}
\overline{B'_2(0)}\subset \bigcup_{j=1}^\tau B_{1}'(x'_j)~\!.
\end{equation}

Take a number $\eps_0$ such that
$0<\eps_0<\frac12 {\mathcal S_\lambda}$.
Notice that the ratio in (\ref{eq:problem}) is invariant with respect to translations in $\R^{n-1}$ and
with respect to the transforms 
$u(x)\mapsto  \alpha u(\beta x)$ for $\alpha\neq 0, \beta>0$. Thus we can
select a bounded minimizing sequence $u_h$ for $\mathcal S_\lambda$ satisfying the normalization condition
\begin{equation}
\label{eq:normaliz}
\|x_1^{-b}u_h\|^{p}_{p}={\mathcal S_\lambda^\frac{p}{p-2}}~,~~
{\mathcal E_s^\lambda}(u_h)={\mathcal S_\lambda^\frac{p}{p-2}}+o(1)
\end{equation}
and such that
\begin{equation}
\label{eq:Rda_sopra}
\eps_0^{\frac{p}{p-2}}\le \max\limits_j\int\limits_0^{2}\int\limits_{B_2'({x'_j})}
x_1^{-pb}|u_h|^{p}~dx'dx_1
\le \int\limits_0^{2}\int\limits_{B'_2(0)}x_1^{-pb}|u_h|^{p}~dx'dx_1
\le   \left(2\eps_0\right)^{\frac{p}{p-2}}.
\end{equation}
Up to a subsequence, we have that $u_h\to u$ weakly in  ${\widetilde{\mathcal D}^s(\R^n_+)}$.
We claim that $u\neq 0$, that is enough to conclude the proof.

Assume by contradiction that
$u=0$. 
By Ekeland's variational principle we can assume that
there exists a sequence ${f}_h\to 0$ in  ${\widetilde{\mathcal D}^s(\R^n_+)}'$,
such that
\begin{equation}
\label{eq:Rtu_equation}
\Ds u_h-\lambda x_1^{-2s}u_h=
x_1^{-pb}|u_h|^{{p}-2}u_h+{f}_h\qquad\textrm{in ~$\widetilde{\mathcal D}^s(\R^n_+)'$}.
\end{equation}

Take a cut-off function $\f\in\mathcal C^\infty_0(-2,2)$ such that 
$\f\equiv 1$ on $(-1,1)$ and define $\f_j(x')=\f(|x'-x'_j|)$, $j=1,\dots,\tau$.

Note that the cut-off function
$\psi_j(x_1,x'):=\f(x_1)\f_j(x')$ has compact support in $(-2,2)\times B'_2(x'_j)$
and that $\psi_j^2u_h$ is a bounded sequence in $\widetilde{\mathcal D}^s(\R^n_+)$ by Lemma \ref{L:technicalD}.
Use $\psi_j^2u_h$ as test function
in (\ref{eq:Rtu_equation}) to find
\begin{equation}
\label{eq:Rmultline}
\langle\Ds u_h,\psi_j^2u_h\rangle-\lambda\irnplus x_1^{-2s}|\psi_j u_h|^2~\!dx
=\int\limits_{\R^n} x_1^{-pb}|u_h|^{{p}-2}|\psi_j u_h|^2~dx+o(1)~\!.
\end{equation}
Thanks to H\"older inequality and (\ref{eq:Rda_sopra}) we can estimate the right-hand side by
\begin{multline}
\int\limits_{\R^n} \!x_1^{-pb}|u_h|^{{p}-2}|\psi_j u_h|^2dx\\
\le \Big(\int\limits_0^{2}\!\!\!\int\limits_{~\!B_2'(x'_j)}\!\!\!x_1^{-pb}|u_h|^{p} dx'dx_1\Big)^{\!\frac{p-2}{p}}
\|x_1^{-b}\psi_ju_h\|^2_{p}\le 2\eps_0 \|x_1^{-b}\psi_ju_h\|^2_{p}~\!.
\label{eq:numero}
\end{multline}
To handle the left-hand side of (\ref{eq:Rmultline}) we use
Lemma \ref{L:technicalD}, the compactness
of embedding $\widetilde{\mathcal D}^s(\R^n_+)\hookrightarrow L^2_{\rm loc}(\R^n)$
and the definition of $\mathcal S_\lambda= 
\mathcal S_s^{\lambda,p}(\R^n_+)$
to obtain
$$
\langle\Ds u_h,\psi_j^2u_h\rangle -\lambda\irnplus x_1^{-2s}|\psi_j u_h|^2~\!dx
={\mathcal E_s^\lambda}(\psi_ju_h)+o(1)
\ge {\mathcal S_\lambda}\|x_1^{-b}\psi_ju_h\|^2_{p}+o(1).
$$
In this way, from (\ref{eq:Rmultline}) we infer
\begin{equation}
{\mathcal S_\lambda}\|x_1^{-b}\psi_ju_h\|^2_{p} \le 2\eps_0 \|x_1^{-b}\psi_ju_h\|^2_{p}+o(1)~\!.
\label{eq:Rchain}
\end{equation}
Since $2\eps_0<{\mathcal S_\lambda}$,  formula (\ref{eq:Rchain}) implies that $\|x_1^{-b}\psi_ju_h\|_{p}=o(1)$. 
But then, using (\ref{eq:Rcover}) and recalling that $\psi_j\equiv 1$ on $(0,1)\times B'_1(x'_j)$, we obtain
\begin{eqnarray*}
\int\limits_0^1\!\!\int\limits_{B'_2(0)}x_1^{-pb}|u_h|^{p}~\!dx'dx_1
&\le& \sum_{j=1}^\tau~\int\limits_0^1\!\!\int\limits_{B'_1(x'_j)}x_1^{-pb}|u_h|^{p}~\!dx'dx_1\\
&\le&
\sum_{j=1}^\tau~
\irnplus x_1^{-pb}|\psi_ju_h|^{p} dx 
=o(1).
\end{eqnarray*}
Comparing with the first inequality in (\ref{eq:Rda_sopra}) we arrive at
\begin{equation}
\label{eq:Rcappa}
2^{-pb}\int\limits_1^2\!\!\int\limits_{B'_2(0)}|u_h|^{p}~dx'dx_1\ge
\int\limits_1^2\!\!\int\limits_{B'_2(0)}x_1^{-pb}|u_h|^{p}~dx'dx_1\ge\eps_0^{\frac{p}{p-2}}+o(1),
\end{equation}
that contradicts the compactness
of embedding $\widetilde{\mathcal D}^s(\R^n_+)\hookrightarrow L^p_{\rm loc}(\R^n)$, as $p<\sstar$.
\QED

\noindent{\bf Proof of Theorem \ref{T:positive}.}
Take any nontrivial function $\f\in\mathcal C^\infty_0(B)$, where $B$ is the unit ball about the origin. Let
$e_1=(1,0,\cdots,0)\in\R^n_+$ and take $h\ge 1$.  Testing
$\mathcal S_s^{\lambda,\sstar}(\R^n_+)$ with $\f_h(x)=\f(h(x-e_1))\in\mathcal C^\infty_0(\R^n_+)$ we obtain
\begin{equation}
\label{eq:BN}
\mathcal S^{\lambda,\sstar}_s(\R^n_+)\le
\frac{\langle\Ds \f,\f\rangle-\lambda\|(x_1+h)^{-s}\f\|^2_{2,B}}
{\|\f\|^2_{\sstar,B}}.
\end{equation}
Letting $h\to \infty$ we infer
$\displaystyle{\mathcal S^{\lambda,\sstar}_s(\R^n_+)\le 
\frac{\langle\Ds \f,\f\rangle}
{\|\f\|^2_{\sstar,B}}}$.
Since $\f$ was arbitrarily chosen we can conclude that
 $$
\mathcal S^{\lambda,\sstar}_s(\R^n_+)\le \inf_{\scriptstyle \f\in\mathcal C^\infty_0(B)\atop\scriptstyle  \f\ne 0}
\frac{\langle\Ds \f,\f\rangle}
{\|\f\|^2_{\sstar,B}}=\mathcal S_s,
$$
and $i)$ is proved.

If $\lambda\le 0$ then trivially $\mathcal S_s^{\lambda,\sstar}(\R^n_+)\ge \mathcal S_s$, because
of (\ref{eq:Sobolev}) holds
for $u\in \widetilde{\mathcal D}^s(\R^n_+)$. Hence 
$\mathcal S_s^{\lambda,\sstar}(\R^n_+)= \mathcal S_s$ and is not attained. 

If $\lambda>0$ we take $h=1$ in (\ref{eq:BN}) to get
$$
\mathcal S^{\lambda,\sstar}_s(\R^n_+)\le \inf_{\scriptstyle \f\in\mathcal C^\infty_0(B)\atop\scriptstyle  \f\ne 0}
\frac{\langle\Ds \f,\f\rangle-2^{-2s}\lambda \|\f\|^2_{2,B}}
{\|\f\|^2_{\sstar,B}}~\!.
$$
Therefore we can use Theorems 4.2, 4.3 in \cite{MN-BN}, see also \cite{SV2}, that give
$\mathcal S_s^{\lambda,\sstar}(\R^n_+)<\mathcal S_s$  if $n\ge 4s$ or if $\lambda>0$ is large enough.
\QED

\noindent 
{\bf Proof of Theorem \ref{T:Dirichlet}.}
The first part of the proof goes as for Theorem
\ref{T:subcritical}.
We assume that $n\ge 2$ and use the same notation as in the proof of Theorem 
\ref{T:subcritical}, with $p=\sstar$ and $b=0$.  

We select a minimizing sequence $u_h$ satisfying (\ref{eq:normaliz}) and (\ref{eq:Rda_sopra}).
Up to a subsequence, we have that $u_h\to u$ weakly in  ${\widetilde{\mathcal D}^s(\R^n_+)}$.
If $u=0$ then we can assume that
there exists a sequence ${f}_h\to 0$ in  ${\widetilde{\mathcal D}^s(\R^n_+)}'$,
such that $u_h$ solves
\begin{equation}
\label{eq:Rtu_equation2}
\Ds u_h-\lambda x_1^{-2s}u_h=
|u_h|^{{\sstar}-2}u_h+{f}_h\qquad\textrm{in ~$\widetilde{\mathcal D}^s(\R^n_+)'$},
\end{equation}
compare with (\ref{eq:Rtu_equation}).
Arguing as in the proof of Theorem \ref{T:subcritical} one can prove that
(\ref{eq:Rcappa}) holds with $p=\sstar$ and $b=0$.

Now we take a cut-off function $\phi\in\mathcal C^\infty_0(\R^n_+)$ such that
$\phi\equiv 1$ on $(1,2)\times B'_2(0)$. We test (\ref{eq:Rtu_equation2}) with $\phi^2u_h\in\widetilde{\mathcal D}^s(\R^n_+)$ to get
\begin{equation}
\label{eq:spring}
\langle\Ds u_h,\phi^2u_h\rangle-\lambda\irnplus x_1^{-2s}|\phi u_h|^2~\!dx
=\int\limits_{\R^n} |u_h|^{\sstar-2}|\phi u_h|^2~dx+o(1)~\!.
\end{equation}
Since $\text{supp}(\phi)\subset \R^n_+$, by compactness
of embedding $\widetilde{\mathcal D}^s(\R^n_+)\hookrightarrow L^2_{\rm loc}(\R^n)$ we have
$\|x_1^{-s}\phi u_h\|_2\to 0$. Thus, we can use Lemma
\ref{L:technicalD} and the Sobolev inequality to infer
$$
\langle\Ds u_h,\phi^2u_h\rangle=\mathcal E_s^\lambda(\phi u_h)+o(1)=
\langle\Ds \phi u_h,\phi u_h\rangle+o(1)
\ge \mathcal S_s\|\phi u_h\|^2_{\sstar}+o(1).
$$
Therefore, estimating the right hand side of (\ref{eq:spring}) via H\"older inequality we obtain
\begin{equation}
\label{eq:miami}
\mathcal S_s\|\phi u_h\|^2_{\sstar}\le 
\|u_h\|^{\sstar-2}_{\sstar}\|\phi u_h\|^2_{\sstar}+o(1)
={\mathcal S_\lambda}\|\phi u_h\|^2_{\sstar}+o(1).
\end{equation}
Now we recall that $S_\lambda<\mathcal S_s$ and
$\phi\equiv 1$ on $(1,2)\times B'_2(0)$. Thus (\ref{eq:miami}) gives  
$$
\int\limits_1^2\!\!\int\limits_{B'_2(0)}|u_h|^{\sstar}~\!dx_1dx'=o(1)~\!.$$ 
We reached a contradiction with (\ref{eq:Rcappa}),
that concludes the proof.
\QED

\section{Additional remarks and  problems}
\label{S:HSM}

In this section we compare the available results for $s\in(0,1)$ with some known results 
in the local case $s=1$, $n\ge 2$, when 
$\mathcal H_1=\frac14$,
$\Ds=-\Delta$ is the standard Laplace operator, and $\langle -\Delta u,u\rangle=\|\nabla u\|_2^2$
for  $u\in {\mathcal D}^1(\R^n)$.

Recall that Maz'ya proved in \cite[2.1.6, Corollary 3]{Ma},
that there exists
a positive best constant $\mathcal S_1^{\frac14,p}(\R^n_+)$ such that 
\begin{equation}
\label{eq:Mazya}
\langle-\Delta u-\mathcal H_1x_1^{-2}u,u\rangle=
\irnplus\big(|\nabla u|^2-\frac14x_1^{-2}|u|^2\big)~\!dx
\ge \mathcal S_1^{\frac14,p}(\R^n_+)
\Big(\irnplus x_1^{-pb}|u|^p~\!dx\Big)^{\frac 2p}
\end{equation}
for any $u\in \mathcal C^\infty_0(\R^n_+)$, where $n\ge3$,
$2<p\le 2^*_1=\frac{2n}{n-2}$ and $\frac{b}{n}=\frac{1}{p}-\frac{n-2}{2n}$, accordingly with (\ref{eq:assumption}).
Inequality
(\ref{eq:Mazya}) holds as well if $n=2$, for any $p>2$ and for $b=\frac 2p$, see \cite[Appendix B]{MS}.

As concerns the attainability of $\mathcal S_1^{\frac14,p}(\R^n_+)$ we refer to
\cite{TT} for $p=2^*_1$ and $n\ge 4$, and to \cite[Sec. 6]{MS} for $p<2^*_1$ and $n\ge2$.
Finally, it was proved in \cite{M} that the best constant $\mathcal S_1^{\lambda,p}(\R^n_+)$  
is attained if $2<p<2^*_1$ and $-\infty<\lambda<\frac 14$,  and when $p=2^*_1$, $n\ge4$ and $0<\lambda<\frac14$
(clearly, $\mathcal S_1^{\lambda,2^*_1}(\R^n_+)$ is never achieved if $\lambda\le 0$).

Surprisingly, in the lower dimensional critical case $n=3$, $p=6$ one has 
$\mathcal S_1^{\lambda,6}(\R^n_+)=\mathcal S_1$ and 
the minimizer never exists, whatever $\lambda\le \frac14$ is (see \cite{BFL} and \cite{MS}).

Now take $s\in(0,1)$, $u\in \mathcal C^\infty_0(\R^n_+)$ and compute
\begin{equation}
\label{eq:F}
\langle \Ds u,u\rangle=\frac{C_{n,s}}{2}\iirnp\frac{(u(x)-u(y))^2}{|x-y|^{n+2s}}~dxdy+\gamma_s\irnplus x_1^{-2s}u^2~\!dx,
\end{equation}
where
$$
\gamma_s=\frac{2^{2s-1}\Gamma\Big(s+\frac12\Big)}{\sqrt\pi~\!\Gamma\big(1-s\big)}.
$$
From the proof of \cite[Lemma 2]{BD} one gets that $\mathcal H_s>\gamma_s$ for $s\neq \frac12$, while $\mathcal H_{\frac12}=\mathcal \gamma_{\frac12}=\frac1\pi$.

The above computation and the inequality proved by C.A. Sloane in 
\cite{Sl} readily imply the next result.

\begin{Proposition}
\label{T:Sl}
Let $n\ge 2$, $s\in(\frac12,1)$. There exists a  best constant $\mathcal S_s^{\mathcal H_s,\sstar}(\R^n_+)>0$ such that
\begin{equation}
\label{eq:Sl}
\langle \Ds u-\mathcal H_sx_1^{-2s}u,u\rangle\ge \mathcal S_s^{\mathcal H_s,\sstar}(\R^n_+)
\Big(\irnplus |u|^\sstar~\!dx\Big)^{\frac 2\sstar}\quad\text{for any  $u\in \mathcal C^\infty_0(\R^n_+)$.}
\end{equation}
\end{Proposition}
In the first version of the present paper the following question has been raised up.
\bigskip

\noindent
{\bf Problem 1.} {\em Let $n\ge 2$, and $p\in(2,\sstar]$. Find sharp conditions on $s\in(0,1)$
that guarantee the existence of a best constant $\mathcal S_s^{\mathcal H_s,p}(\R^n_+)>0$ such that for
$b=b(n,s,p)$ as in (\ref{eq:assumption}) one has}
$$
\langle \Ds u-\mathcal H_sx_1^{-2s}u,u\rangle
\ge \mathcal S_s^{\mathcal H_s,p}(\R^n_+)
\Big(\irnplus x_1^{-pb}|u|^p~\!dx\Big)^{\frac 2p}
\quad\textit{for any $u\in \mathcal C^\infty_0(\R^n_+)$}.
$$
In the recent publication \cite{DLV},  Dyda, Lehrb\"ack and V\"ah\"akangas
gave a complete answer to Problem 1. As far as we know, the next problem is still open. 

\medskip

\noindent
{\bf Problem 2} {\em Assume $\mathcal S_s^{\mathcal H_s,p}(\R^n_+)>0$. Is $\mathcal S_s^{\mathcal H_s,p}(\R^n_+)$ attained?}
\medskip

Inspired by the result of \cite{MS}, we formulate the following conjecture.
\bigskip

\noindent
{\bf Conjecture} {\em Let $s\in (0,1)$, $2s<n< 4s$ (hence, $n\le 3$). Then the best constant
$\mathcal S_s^{\lambda,\sstar}(\R^n_+)$ is never achieved.}

\footnotesize
\label{References}

\end{document}